\documentclass[a4paper, 11pt]{amsart}
\usepackage{mathptmx,amssymb,amscd,latexsym, eulervm}
\usepackage{amsmath}
\usepackage{amsthm}
\usepackage{mathdots}
\usepackage[colorlinks=true,citecolor=violet,linkcolor=blue,urlcolor=blue]{hyperref}
\usepackage[dvipsnames]{xcolor}
\usepackage[onehalfspacing]{setspace}
\usepackage{tabularx}
\usepackage{amsfonts}
\usepackage{paralist}
\usepackage{aliascnt}
\usepackage{amscd}
\usepackage{blkarray}
\usepackage{mathbbol}
\usepackage{setspace}
\usepackage{needspace}
\usepackage[inner=2.4cm,outer=2.4cm, bottom=3.2cm]{geometry}
\usepackage{tikz, tikz-cd}
\usepackage{calligra,mathrsfs}

\usepackage{tikz}
\usetikzlibrary{matrix}
\usetikzlibrary{arrows,calc}
\allowdisplaybreaks

\AtBeginDocument{%
	\def\MR#1{}
}

\makeatletter
\@namedef{subjclassname@2020}{%
	\textup{2020} Mathematics Subject Classification}
\makeatother

\newcommand{\Bl}{{\mathcal{B}\ell}}

\newcommand{\kk}{\mathbb{k}}

\newcommand{\ZZ}{\mathbb{Z}}

\newcommand{\mm}{{\normalfont\mathfrak{m}}}

\newcommand{\pp}{{\normalfont\mathfrak{p}}}

\newcommand{\qqq}{\normalfont\mathfrak{q}}

\newcommand{\qq}{{\mathfrak{q}}}

\newcommand{\Ker}{\normalfont\text{Ker}}

\newcommand{\HT}{\normalfont\text{ht}}

\newcommand{\Supp}{\normalfont\text{Supp}}
\newcommand{\Ass}{{\normalfont\text{Ass}}}

\newcommand{\Rees}{\mathscr{R}}

\newcommand{\OO}{\mathscr{O}}

\newcommand{\LL}{\mathscr{L}}

\newcommand{\HL}{\normalfont\text{H}_{\mm}}
\newcommand{\HH}{\normalfont\text{H}}

\newcommand{\gr}{{\normalfont\text{gr}}}

\newcommand{\Proj}{\normalfont\text{Proj}}

\newcommand{\Spec}{\normalfont\text{Spec}}

\newcommand{\biProj}{{\normalfont\text{BiProj}}}

\newcommand{\wt}{\widetilde}

\def\f0{\mathbf{0}}

\def\1{\mathbf{1}}

\newtheorem{theorem}{Theorem}[section]

\newtheorem{headthm}{Theorem}

\newaliascnt{headcor}{headthm}

\aliascntresetthe{headcor}

\newaliascnt{headconj}{headthm}

\aliascntresetthe{headconj}

\newaliascnt{corollary}{theorem}
\newtheorem{corollary}[corollary]{Corollary}
\aliascntresetthe{corollary}

\newaliascnt{claim}{theorem}

\aliascntresetthe{claim}

\newaliascnt{lemma}{theorem}
\newtheorem{lemma}[lemma]{Lemma}
\aliascntresetthe{lemma}

\newaliascnt{conjecture}{theorem}

\aliascntresetthe{conjecture}

\newaliascnt{proposition}{theorem}
\newtheorem{proposition}[proposition]{Proposition}
\aliascntresetthe{proposition}

\theoremstyle{definition}
\newaliascnt{definition}{theorem}
\newtheorem{definition}[definition]{Definition}
\aliascntresetthe{definition}

\newaliascnt{notation}{theorem}

\aliascntresetthe{notation}

\newaliascnt{convention}{theorem}
\newtheorem{convention}[convention]{Convention}
\aliascntresetthe{convention}

\newaliascnt{example}{theorem}
\newtheorem{example}[example]{Example}
\aliascntresetthe{example}

\newaliascnt{examples}{theorem}

\aliascntresetthe{examples}

\newaliascnt{remark}{theorem}
\newtheorem{remark}[remark]{Remark}
\aliascntresetthe{remark}

\newaliascnt{question}{theorem}

\aliascntresetthe{question}

\newaliascnt{questions}{theorem}

\aliascntresetthe{questions}

\newaliascnt{problem}{theorem}

\aliascntresetthe{problem}

\newaliascnt{construction}{theorem}

\aliascntresetthe{construction}

\newaliascnt{setup}{theorem}
\newtheorem{setup}[setup]{Setup}
\aliascntresetthe{setup}

\newaliascnt{algorithm}{theorem}

\aliascntresetthe{algorithm}

\newaliascnt{observation}{theorem}

\aliascntresetthe{observation}

\newaliascnt{defprop}{theorem}

\aliascntresetthe{defprop}

\newaliascnt{fact}{theorem}

\aliascntresetthe{fact}

\DeclareFontFamily{OT1}{pzc}{}
\DeclareFontShape{OT1}{pzc}{m}{it}{<-> s * [1.100] pzcmi7t}{}
\DeclareMathAlphabet{\mathchanc}{OT1}{pzc}{m}{it}

\def\equationautorefname~#1\null{(#1)\null}
\def\sectionautorefname~#1\null{Section #1\null}
\def\subsectionautorefname~#1\null{\S #1\null}

\def\trdeg{{\rm trdeg}}
\def\surjects{\twoheadrightarrow}

\title{Intersection Product Formulas and Integral Dependence}
\author{Yairon Cid-Ruiz}
\address[Cid-Ruiz]{Department of Mathematics, North Carolina State University, Raleigh, NC 27695, USA}
\email{ycidrui@ncsu.edu}

\author{Steven Dale Cutkosky}
\address[Cutkosky]{Department of Mathematics, University of Missouri, Columbia, MO 65211, USA}
\email{cutkoskys@missouri.edu}

\date{\today}

\begin{document}
	\begin{abstract}
		We prove the intersection product formulas of Gaffney and Gassler for polar multiplicities and Segre numbers over arbitrary Noetherian local rings. 
		For equidimensional and universally catenary rings, we give a geometric refinement of the Principle of Specialization of Integral Dependence, in which the maximal dimension of the images of exceptional components mapping into the hypersurface determines a precise Segre number that strictly increases.
	\end{abstract}
	
	\maketitle
	
	\vspace*{-.6cm}
	
	\section{Introduction}
	\label{sec:introduction}
	
	The study of integral dependence through numerical invariants goes back to Rees' multiplicity theorem \cite{Rees61}.
	In an equidimensional and universally catenary Noetherian local ring $(R,\mm)$, this theorem states that two $\mm$-primary ideals $I\subseteq J$ have the same integral closure if and only if $e(I,R)=e(J,R)$. 
	For ideals of arbitrary height, the corresponding criterion is expressed in terms of Segre numbers. 
	In the complex analytic setting, Gaffney and Gassler \cite{GG} proved this extension of Rees' theorem by establishing a Principle of Specialization of Integral Dependence. 
	Their work also establishes the lexicographic upper semicontinuity of Segre numbers and gives cycle descriptions and intersection product formulas for Segre numbers and polar multiplicities.
	
	Several of Gaffney--Gassler's results have since been extended to the algebraic setting.
	Achilles and Rams \cite{AR} identified the Segre numbers of Gaffney and Gassler with the multiplicity sequence introduced by Achilles and Manaresi \cite{AM}.
	Polini, Trung, Ulrich, and Validashti \cite{PTUV} proved that equality of these sequences characterizes integral dependence of arbitrary ideals in an equidimensional and universally catenary Noetherian local ring.
	Algebraic cycle descriptions and semicontinuity results were obtained by Cid-Ruiz, Polini, and Ulrich \cite{CRPU}. 
	One of the main goals of this paper is to establish the intersection product formulas of Gaffney and Gassler over arbitrary Noetherian local rings. 
	We also use the geometry of the blow-up to give a different proof of the specialization criterion in \cite[Theorem~3.3]{PTUV} and obtain a more precise statement that identifies a precise Segre number that must increase.
	
	\smallskip
	
	We now describe our results more precisely. 
	Let $(R,\mm,\kappa)$ be a Noetherian local ring of dimension $d$ and $I\subset R$ be a proper ideal.
	Let $X = \Spec(R)$ and consider the punctured spectrum $U = X \setminus \{\mm\}$.
	The \emph{polar multiplicities} $m_i(I,R)$ and the \emph{Segre numbers} $c_i(I,R)$ are defined by the bigraded Hilbert polynomials associated to the standard bigraded $\kappa$-algebras
	$$
	\gr_\mm\left(\Rees(I)\right)
	\qquad\text{and}\qquad
	\gr_\mm\left(\gr_I(R)\right),
	$$
	respectively (see \autoref{def_polar_segre}). 
	We also consider the \emph{polar-Segre multiplicities} $\nu_i(I,R)=m_i(I,R)+c_i(I,R)$.
	
	To express these numbers as intersection products, we work on the joint
	blow-up
	$$
	b\colon \; B=\Bl_{\mm I}(X) \;\longrightarrow\; X.
	$$
	On $B$, both $I\OO_B$ and $\mm\OO_B$ are invertible ideals, and so they determine effective Cartier divisors $E_I$ and $E_\mm$ given by $I\OO_B=\OO_B(-E_I)$ and $\mm\OO_B=\OO_B(-E_\mm)$.
	We write
	$
	\left[E_I\right] =  \left[E_I\right]^{b^{-1}(U)}+\left[E_I\right]^{b^{-1}(\{\mm\})},
	$
	where the two summands consist of the components meeting $b^{-1}(U)$
	and those mapping to the closed point, respectively, with their multiplicities.
	The following theorem expresses the multiplicities in terms of these
	divisors and cycles.
	
	\begin{headthm}[{Gaffney-Gassler intersection formulas;~\autoref{thm:main}}]
		\label{head:intersection-formulas}
		With the notation above, for all $1 \le i \le d-1$, the following formulas hold:
		\begin{align*}
			m_i(I,R)&\;=\;-\int_B (-E_I)^i(-E_\mm)^{d-i}\cdot[B],\\
			c_i(I,R)&\;=\;-\int_B (-E_I)^{i-1}(-E_\mm)^{d-i}\cdot[E_I]^{b^{-1}(U)},\\
			\nu_i(I,R)&\;=\;\int_B (-E_I)^{i-1}(-E_\mm)^{d-i}\cdot[E_I]^{b^{-1}(\{\mm\})}.
		\end{align*}
		Moreover, when $\kappa$ is infinite, for all $1 \le i \le d-2$, we have the following cutting formulas: 
		\begin{enumerate}[\rm (i)]
			\item If $\dim(R/I) < d$ and $f \in I$ is a general element, then 
			$$
			m_i(I,R/fR) \;=\; m_{i+1}(I,R), \quad\;\;
			c_i(I,R/fR) \;=\; c_{i+1}(I,R), \quad\;\;
			\nu_i(I,R/fR) \;=\; \nu_{i+1}(I,R).
			$$
			\item If $x \in \mm$ is a general element, then 
			$$
			m_i(I,R/xR) \;=\; m_{i}(I,R), \quad\; c_i(I,R/xR) \;=\; c_{i}(I,R), \quad\; \nu_i(I,R/xR) \;=\; \nu_{i}(I,R).
			$$
		\end{enumerate}
	\end{headthm}
	
	Here the integrals denote local intersection degrees with support over the closed point (see \autoref{subsec:thorup-intersection}). 
	The proof of \autoref{head:intersection-formulas} uses Thorup's rational equivalence theory \cite{Thorup}, which allows us to work over an arbitrary Noetherian local ring.

	\smallskip
	
We next turn to specialization, the main step in the integral-dependence criteria established in \cite{GG} and \cite{PTUV}.
	Geometrically, the issue is	whether a component of the exceptional divisor $\Proj(\gr_I(R))$ maps into a given hypersurface $V(t)\subset X$.
	Our next result shows that such components force a strict increase in a precise Segre number: the index of this number is determined by the largest dimension of their images in $X$.
	
	\begin{headthm}[\autoref{thm:avoidance_main_explicit}]
		\label{head:strict-specialization}
		Assume that $R$ is equidimensional and universally catenary and that $\HT(I) \ge 1$.
		Let $t\in\mm$ be an element such that $\HT(t) = 1$ and $\HT(I, t) \ge 2$.
		Assume there is a minimal prime of $\gr_I(R)$ that contains $t$. 
		Let $s$ be the maximal dimension
		$$
		s \;:=\; \max\big\lbrace \dim(R/Q \cap R) \mid \text{$Q$ minimal prime of $\gr_I(R)$ with $t \in Q$}\big\rbrace
		$$
		of the contractions to $R$ of the minimal primes of $\gr_I(R)$ containing $t$.
		Let $j := d-1-s$.
		Then we have the strict inequality
		$$
		c_j(I, R/tR) \;>\; c_j(I, R).
		$$
	\end{headthm}

	In an analytic setting, Teissier’s Principle of Specialization of Integral Dependence (PSID) (see \cite{TEISSIER_CYC} and \cite[Appendice~I]{TEISSIER_RES2}) states that for a family of zero-dimensional ideals with constant Hilbert-Samuel multiplicity, a section is integrally	dependent on the total family if and only if it is integrally dependent on the fibers corresponding to a Zariski-open dense subset of the base. 
	Using Segre numbers, Gaffney and Gassler \cite{GG} (also in the analytic setting) extended the PSID to families of ideals that are not necessarily zero-dimensional.	
	As a direct consequence of \autoref{head:strict-specialization}, we obtain in \autoref{cor_PSID} an algebraic version of the PSID, showing that, under certain hypotheses, it suffices to check integral dependence on the open set $D(t)\subset \Spec(R)$.

	\smallskip
	
	The basic outline of this paper is as follows.
	
	\medskip
	\noindent
	\textbf{Outline.}
	In \autoref{sec:preliminaries}, we recall the definitions and preliminary
	results used throughout the paper. The proof of
	\autoref{head:intersection-formulas}, together with the cutting formulas
	and cycle descriptions, is given in \autoref{sec:intersection-formulas}.
	We conclude that section by studying integral dependence under general hypersurface sections (see \autoref{subsec:integral-dependence-general-sections}).
	Finally, in \autoref{sec:integral-dependence}, we prove
	\autoref{head:strict-specialization}.

	\section{Basic tools}
	\label{sec:preliminaries}
	\label{sec:intersection-general-sections}
	In this section, we recall the definitions and results used throughout the
	paper. 
	We also prove some basic tools that will be necessary.
	The following data is used throughout.

	\begin{setup}
		\label{setup_basic}
		Let $(R,\mm,\kappa)$ be a $d$-dimensional Noetherian local ring and
		$I\subset R$ be a proper ideal.
		Let $X:=\Spec(R)$.		
	\end{setup}
	
We denote by $Z_i(X)$ the free Abelian group of cycles on $X$ of Krull dimension $i$ (in \autoref{subsec:thorup-intersection} below, we also consider the notion of rational equivalence developed by Thorup~\cite{Thorup}).
	For $r\geq0$, set $Z_{\leq r}(X):=\bigoplus_{i=0}^r Z_i(X)$.
	We define the multiplicity homomorphism $e_r\colon Z_{\leq r}(X)\to\ZZ$ by
	$$
	e_r\left(\sum_{\pp}n_{\pp}\left[V(\pp)\right]\right)
	\;:=\;\sum_{\dim(R/\pp)=r}n_{\pp}\,e(R/\pp),
	$$
	where $e(R/\pp)$ denotes Hilbert--Samuel multiplicity.

	\subsection{Polar multiplicities and Segre numbers}
	First, we recall the notions of polar multiplicities and Segre numbers. 	
	
	\begin{definition}
		\label{def_polar_segre}
		Consider the standard bigraded $\kappa$-algebras
		$
		A :=\gr_\mm(\Rees(I))
		$
		and
		$
		G := \gr_\mm(\gr_I(R)).
		$
		Their bigraded components are
		$$
		A_{v,n} \;=\; \frac{\mm^nI^v}{\mm^{n+1}I^v} \qquad \text{and} \qquad
		G_{v,n} \;=\; \frac{\mm^nI^v+I^{v+1}}{\mm^{n+1}I^v+I^{v+1}}.
		$$
		For $v,n\gg0$, the Hilbert function of $A$ is a polynomial with homogeneous part of degree $d-1$ given by
		$$
		\sum_{i=0}^{d-1}m_i(I,R)\frac{v^in^{d-1-i}}{i!(d-1-i)!}.
		$$
		The second sum transform of the Hilbert function
		$$
		H_G(v,n) \;:=\; \sum_{a=0}^{v}\sum_{b=0}^{n}\dim_\kappa\left(G_{a,b}\right)
		$$
		is eventually a polynomial with homogeneous part of degree $d$ given by
		$$
		\sum_{i=0}^{d}c_i(I,R)\frac{v^in^{d-i}}{i!(d-i)!}.
		$$
		We call $m_i(I,R)$ the \emph{polar multiplicities} and $c_i(I,R)$ the
		\emph{Segre numbers} of $I$.
		We put $m_d(I,R)=0$ and
		$\nu_i(I,R)=m_i(I,R)+c_i(I,R)$ for $1\leq i\leq d$.
	\end{definition}
	
	When $d\geq2$, the coefficients $c_i(I,R)$ with $1\leq i\leq d-1$ can be read from the ordinary bigraded Hilbert polynomial of $G$.
	The two end coefficients are given by
	$$
	c_0(I,R) \;=\; e(R)-m_0(I,R)  \qquad \text{and} \qquad
	c_d(I,R)\;=\; j(I)
	\;=\; e_{d-1}\bigl(\HL^0(\gr_I(R))\bigr).
	$$
	We have that: if $\HT(I)\geq1$, then $c_0(I,R)=0$;
	if $I$ is $\mm$-primary, then $c_i(I,R)=0$ for $i<d$ and
	$c_d(I,R)=e(I,R)$.

	\subsection{Thorup's rational equivalence theory}
	\label{subsec:thorup-intersection}
	
	We recall the basics of the rational equivalence theory developed by Thorup \cite{Thorup} over Noetherian schemes, which will be used in this paper (a slightly less general treatment is made in \cite[\href{https://stacks.math.columbia.edu/tag/02QK}{Tag 02QK}]{stacks-project} over a universally catenary scheme).  
	See also the recent developments \cite{CutMon, CutMon2,  Cut26}.
	
	A \emph{grading} on a Noetherian base scheme assigns integers to its integral
	closed subschemes and induces a grading on schemes essentially of finite type over the base.  
	In Thorup's theory \cite{Thorup} one can always use the \emph{canonical grading} given by
	$$
	\delta_{\rm can}(V) \;:=\; -\dim\left(\OO_{X,\eta_V}\right),
	$$
	where $\eta_V$ is the generic point of an integral closed subscheme $V \subset X$.  
	Since $X=\Spec(R)$ has finite Krull dimension, we use instead the \emph{topological
		grading}
	$$
	\delta_{\rm top}(V) \;:=\; \dim\left(V\right)
	$$
	for an integral closed subscheme $V \subset X$.

	Let $\pi \colon Y\to X$ be a morphism essentially of finite type.
	The grading on $Y$ \emph{induced} by the topological grading on $X$ is given by
	$$
	\delta_{X, \rm top}(W) \;:=\; \trdeg_{\kappa(V)}\left(\kappa(W)\right) +\dim\left(V\right),
	$$
	where $W\subset Y$ is
	an integral closed subscheme and $V=\overline{\pi(W)} \subset X$.  
	The dimension formula gives 
	$$
	\dim(W) \;\leq\; \delta_{X, \rm top}(W)
	$$ 
	(see, e.g.,
	\cite[Theorem~15.5]{MatCRT}). 
	Moreover, Thorup \cite[Proposition~4.3]{Thorup} showed that
	$$
	\text{if $\pi$ is proper, \qquad then} \quad  \delta_{X, \rm top}(W) \;=\; \dim(W).
	$$

	Given a grading on $X$ and a scheme $Y$ essentially of finite type over $X$, by utilizing the \emph{induced grading} on $Y$, Thorup builds a rational equivalence theory for $Y$.
	This leads to the following very convenient  convention. 
	
	\begin{convention}
		Let $Y$ be a proper scheme over $X = \Spec(R)$.
		Throughout this paper, we consider Thorup's rational equivalence theory on $Y$ with the induced grading $\delta_{X, \rm top}$.
		Most importantly, the topological grading of $Y$ coincides with the induced grading $\delta_{X, \rm top}$ on $Y$.
	\end{convention}

	Therefore, for every proper $X$-scheme $Y$, we may use the ordinary dimensional cycle group $Z_r(Y)$.  
	We write $A_r(Y/X)$ for the quotient by Thorup rational equivalence and $[\alpha]_r$ for the $r$-dimensional part of a cycle $\alpha$.  
	If $f : Z \rightarrow Y$ is a proper morphism of proper $X$-schemes, we get a functorial homomorphism
	$$
	f_*\colon \; A_r(Z/X) \;\longrightarrow\; A_r(Y/X)
	$$
	by \cite[Proposition~6.5]{Thorup}. 
	Moreover, Cartier-divisor intersections and the	projection formula are available for all the proper schemes used	throughout this paper.
	
	We use the following convention for local intersection degrees.
	Let $\pi\colon Y\to X$ be a proper morphism and set $Y_0:=Y\times_X\Spec(\kappa)$.
	Since $Y_0$ is proper over $\Spec(\kappa)$, its structure morphism $\pi_0\colon Y_0\to\Spec(\kappa)$ induces the usual degree homomorphism
	$$
	\deg_\kappa\colon\; A_0\left(Y_0\right) \;\xrightarrow{\;(\pi_0)_*\;}\; A_0\left(\Spec(\kappa)\right) \;\cong\; \ZZ
	$$
	(see \cite[Definition~1.4]{Fu}).
	For a class $\beta\in A_0(Y_0)$ represented by a zero-dimensional cycle $\sum_V n_V[V]$, we write
	$$
	\int_Y\beta \;:=\; \deg_\kappa(\beta)
	\;=\; \sum_V n_V\left[R(V):\kappa\right].
	$$
	This integer is independent of the chosen representative, since proper pushforward respects rational equivalence on $Y_0$.

	\subsection{Passage to an infinite residue field}
	\label{subsec:infinite-residue-extension}
	
	Let $S:=R[z]_{\mm R[z]}$, where $z$ is a new indeterminate.
	The extension $R\to S$ is faithfully flat and
	local, $\dim(S)=d$, its maximal ideal is $\mm S$, and its residue field is $\kappa(z)$. 
	For every ideal $I \subset R$, flat base change gives
	$$
	(R/I) \otimes_R S \; \cong \; (R/I)[z]_{\mm(R/I)[z]} 
	\qquad \text{ and } \qquad
	\Bl_{IS}(\Spec(S)) \;\cong \; \Bl_I(X)\times_X\Spec(S).
	$$
	The bigraded components in \autoref{def_polar_segre} are extended from $\kappa$ to $\kappa(z)$, so their dimensions and the resulting multiplicities are unchanged. 
	For a proper scheme $Y$ over $X = \Spec(R)$, on the closed fiber $Y_0 = Y \times_{X} \Spec(\kappa)$, the base change is the field extension
	$\kappa\subset\kappa(z)$, and so it preserves degrees and commutes with proper pushforward and Cartier intersections.
	Therefore, passing to an infinite residue field is a mild operation that keeps the intersection-theoretic data.	
	
	\subsection{General sections}
	
	We first give the avoidance argument used to choose general sections.
	
	\begin{lemma}
		\label{lem:general-section-avoidance}
		Assume that the residue field $\kappa$ is infinite.  
		Let $\pi \colon Y \to X=\Spec(R)$ be a proper morphism and $\mathscr{L}$ be a globally generated line bundle on $Y$.
		Choose points $\eta_1,\ldots,\eta_s \in Y$.
		Then, for a general section $\sigma \in \HH^0(Y, \LL)$, the zero scheme $Z(\sigma) \subset Y$ does not contain any of the points $\eta_1,\ldots,\eta_s$.
	\end{lemma}
	\begin{proof}
		Since $\pi$ is proper and $R$ is Noetherian, $M:=\HH^0(Y,\LL)$ is a finitely generated $R$-module.
		For $1\leq i\leq s$, let
		$$
		M_i \;:=\; \Ker\left(M\longrightarrow
		\LL_{\eta_i}\otimes_{\OO_{Y,\eta_i}}\kappa(\eta_i)\right).
		$$
		The global generation of $\LL$ implies that each $M_i$ is a proper
		submodule of $M$.  Let
		$$
		V \; :=\; M/\mm M
		\qquad \text{ and } \qquad
		V_i \;:=\; (M_i+\mm M)/\mm M\subset V.
		$$
		If $V_i=V$, then Nakayama's lemma applied to $M/M_i$ would give
		$M=M_i$, a contradiction.  Thus each $V_i$ is a proper
		$\kappa$-linear subspace of $V$.  Since $\kappa$ is infinite,
		$V\setminus\bigcup_iV_i$ is a nonempty Zariski-open subset.
		Every section $\sigma\in M$ whose image lies in this open subset has nonzero value in the fiber of $\LL$ at each $\eta_i$, and so $\eta_i\notin Z(\sigma)$ for all $i$.
	\end{proof}
	
	\begin{corollary}
		\label{lem_avoid_projection}
		Assume that the residue field $\kappa$ is infinite.
		Let $\pi\colon Y\to X=\Spec(R)$ be a proper morphism with $Y$ integral and $\mathscr{L}$ be a globally generated line bundle on $Y$.
		Let $T\subset X$ be a closed subset such that $\pi(Y)\not\subset T$.
		Then, for a general section $\sigma\in\HH^0(Y,\LL)$, no irreducible
		component of $\pi\left(Z(\sigma)\right)$ is contained in $T$.
	\end{corollary}
	\begin{proof}
		Let $W:=\pi^{-1}(T)\subsetneq Y$.
		By \autoref{lem:general-section-avoidance}, the zero scheme $Z(\sigma)$
		of a general section $\sigma \in \HH^0(Y, \LL)$ avoids the generic point of $Y$ and the generic points of all irreducible components of $W$.
		Since $Y$ is integral, every irreducible component of $Z(\sigma)$ has codimension one in $Y$ by the principal ideal theorem.
		If such a component were contained in $W$, it would therefore be an irreducible component of $W$, a contradiction.
		So the proof of the corollary is complete.
	\end{proof}
	
	The following lemma describes specialization of a biprojective Rees algebra
	by a general element.
	
	\begin{lemma}
		\label{lem:rees-superficial-specialization}
		Assume that the residue field $\kappa$ is infinite.
		Let $I,J\subseteq R$ be ideals, and consider the multi-Rees algebra $\Rees_R(I,J)=R[IT_1,JT_2]$.
		For a general element $x\in I$, put $\overline R=R/(x)$,
		$\overline I=I\overline R$, and $\overline J=J\overline R$.
		The natural surjection
		$\Rees_R(I,J)/(xT_1)\longrightarrow\Rees_{\overline R}(\overline I,\overline J)$
		induces a natural isomorphism
		$$
		\biProj\bigl(\Rees_R(I,J)/(xT_1)\bigr)
		\;\cong\;
		\biProj\bigl(\Rees_{\overline R}(\overline I,\overline J)\bigr).
		$$
	\end{lemma}
	\begin{proof}
		This follows from Rees' lemma \cite[Lemma~1.2]{Rees84}; see also
		Trung--Verma \cite{TV}.
	\end{proof}

	\section{Intersection product formulas}
	\label{sec:intersection-formulas}

	In this section, by building on the work of Thorup~\cite{Thorup} (see \autoref{subsec:thorup-intersection}), we prove the intersection product formulas stated in \autoref{head:intersection-formulas}.
	Throughout this section, we use the following setup.	
	
	\begin{setup}
		\label{setup_general}
		Let $(R,\mm,\kappa)$ be a $d$-dimensional Noetherian local ring and
		$I\subset R$ be a proper ideal.
		Let $X:=\Spec(R)$ and consider the punctured spectrum $U := X \setminus \{\mm\}$.
		Let
		$
		b_I\colon B_I:=\Bl_I(X)\rightarrow X
		$
		be the blow-up of $X=\Spec(R)$ along $I$. 
		Likewise, the blow-up of $X$ at the origin is
		$
		b_\mm\colon B_\mm:=\Bl_{\mm}(X)\rightarrow X.
		$
		Finally, we consider the joint blow-up
		$
		b\colon B:=\Bl_{\mm I}(X)\rightarrow X.
		$
		Algebraically, we have the descriptions 
		$$
		B_I \;=\; \Proj\big(\Rees(I)\big), \quad  B_\mm \;=\; \Proj\big(\Rees(\mm)\big) \quad \text{ and } \quad B \;=\; \biProj\big(\Rees(I, \mm)\big) \;\cong\; \Proj\big(\Rees(\mm I)\big).
		$$
		Here $\Rees(I) := R[IT] \subset R[T]$, $\Rees(\mm) := R[\mm T] \subset R[T]$ and $\Rees(\mm I) := R[\mm IT] \subset R[T]$ are Rees algebras, and $\Rees(I, \mm) := R[IT_1, \mm T_2] \subset R[T_1, T_2]$ is a multi-Rees algebra. 
		We have the following commutative diagram
		$$
		\begin{tikzcd}[column sep=5.5em, row sep=3.5em]
			B \arrow[r, "\wt b_\mm"] \arrow[d, "\wt b_I"'] \arrow[rd, "b"']
			& B_I \arrow[d, "b_I"] \\
			B_\mm \arrow[r, "b_\mm"']
			& X
		\end{tikzcd}
		$$
		with
		$
		b=b_I\circ \wt b_\mm=b_\mm\circ \wt b_I.
		$
		Let $D_I\subset B_I$ and $D_\mm\subset B_\mm$ be the exceptional divisors of $b_I$ and $b_\mm$, respectively. 
		On the
		joint blow-up $B$, consider the effective Cartier divisors
		$$
		E_I \;:=\; V\left(I \OO_B\right),\qquad E_\mm \;:=\; V\left(\mm \OO_B\right) \qquad \text{and}\qquad E \;:=\; V\left(\mm I \OO_B\right).
		$$
		Then $E$ is the exceptional divisor of $B$, and $E_I=\wt{b}_\mm^*(D_I)$ and $E_\mm=\wt{b}_I^*(D_\mm)$.
		Let
		$$
		h_I \;:=\; c_1(\mathscr O_{B}(1, 0)),
		\qquad
		h_\mm \;:=\; c_1(\mathscr O_{B}(0, 1)),
		\qquad
		h \;:=\; c_1(\mathscr O_B(1,1)).
		$$
	\end{setup}	
	
	We start with the following basic lemma. 
	
	\begin{lemma}
		\label{prop:basic}
		On $B=\Bl_{\mm I}(X)$, the following statements hold:
		\begin{enumerate}[\rm(i)]
			\item $\OO_B(1,0) = \OO_B(-E_I)$, $\OO_B(0,1) = \OO_B(-E_\mm)$ and $\OO_B(1,1) = \OO_B(-E)$.
			\item $E = E_I + E_\mm$.
			\item $h = h_I + h_\mm$.
		\end{enumerate}
	\end{lemma}
	\begin{proof}
		Since $B=\biProj(\Rees(I,\mm))$, the invertible sheaves
		$\OO_B(1,0)$, $\OO_B(0,1)$, and $\OO_B(1,1)$ correspond to the invertible
		ideals $I\OO_B$, $\mm\OO_B$, and $(\mm I)\OO_B$, respectively. This gives
		$$
		\OO_B(1,0)=\OO_B(-E_I), \qquad \OO_B(0,1)=\OO_B(-E_\mm), \qquad
		\OO_B(1,1)=\OO_B(-E).
		$$
		Since $(\mm I)\OO_B=\mm\OO_B\cdot I\OO_B$, we get $E=E_I+E_\mm$. Taking first
		Chern classes gives $h=h_I+h_\mm$.
	\end{proof}
	
	\begin{remark}
		We recall some useful notation from \cite{KTBR,KTMixed}. Let $Y$ be a scheme, $W$ a subscheme of $Y$, and $\mathbf{S}\in Z_*(Y)$ a cycle on $Y$. We define the $W$-part $\mathbf{S}^W$ as follows: if $\mathbf{S}$ is the fundamental cycle of an integral closed subscheme $S\subset Y$, then $\mathbf{S}^W=\mathbf{S}$ if $W$ contains the generic point of $S$, and $\mathbf{S}^W=0$ otherwise; then extend this definition by linearity. 
	\end{remark}

	The next proposition will play an important role in our approach. 
	
	\begin{proposition}
		\label{prop:joint-models}
		On $B = \Bl_{\mm I}(X)$, the following equalities of cycles hold:
		$$
		\left[E_\mm\right] \;=\; \left[\biProj\big(\gr_\mm(\Rees(I))\big)\right]
		\qquad\text{and}\qquad
		\left[E_I\right]^{\,b^{-1}(U)} \;=\; \left[\biProj\big(\Rees(\mm \cdot\gr_I(R))\big)\right].
		$$
	\end{proposition}
	\begin{proof}
		For the first equality, since
		$
		E_\mm = V(\mm\OO_B),
		$
		we have
		$
		\gr_\mm(\Rees(I)) \cong \Rees(I,\mm)/\mm\Rees(I,\mm),
		$
		and therefore
		$
		E_\mm \cong \biProj\big(\gr_\mm(\Rees(I))\big).
		$
		
		For the second equality, since
		$
		E_I = V(I\OO_B),
		$
		we have $E_I\cong \biProj(S) \subset B$ with $	S:=\Rees(I,\mm)/I\Rees(I,\mm)$. 
		Consider the specialization map 
		$$
		\mathfrak{s} \colon \;  S \;\cong\; \frac{\Rees(\mm \cdot \Rees(I))}{I \cdot \Rees(\mm \cdot \Rees(I))} \;\;\surjects\;\; \Rees(\mm \cdot\gr_I(R)) \;\cong\; \Rees\left(\mm \cdot\frac{\Rees(I)}{I\cdot \Rees(I)}\right).
		$$
		For any prime $\pp \in U$ on the punctured spectrum, we have that the map $\mathfrak{s} \otimes_R R_\pp$ is an isomorphism. 
		Hence we obtain 
		$$
		\mm^l \cdot \Ker(\mathfrak{s}) \;=\; 0 \quad \text{ for some $l \ge 1$}.
		$$
		This yields the following set-theoretical decomposition 
		\begin{equation}
			\label{eq_decomp}
			\Supp_{E_I}\left(\widetilde{S/\mm \cdot S}\right) \;\cup\; \Supp_{E_I}\left(\widetilde{S/\Ker(\mathfrak{s})}\right) \;=\; E_I.
		\end{equation}
		Since the biprojective scheme $\biProj\big(\Rees(\mm \cdot\gr_I(R))\big) \subset B$ is isomorphic to the blow-up $\Bl_Z(Y)$ of $Y := \Proj\left(\gr_I(R)\right)$ along $Z := V(\mm \cdot \OO_Y)$, its generic points are contained in $b^{-1}(U)$.
		On the other hand, by the decomposition \autoref{eq_decomp}, the generic points of $E_I$ that are contained in $b^{-1}(U)$ should also be contained in $V(\Ker(\mathfrak{s})) \cong \biProj\big(\Rees(\mm \cdot\gr_I(R))\big)$.
		Therefore, the asserted equality of cycles
		$$
		[E_I]^{\,b^{-1}(U)} \;=\; \left[\biProj\big(\Rees(\mm \cdot\gr_I(R))\big)\right]
		$$
		follows. 
	\end{proof}

	We can now express polar multiplicities and Segre numbers as intersection products of divisors on the joint blow-up $B = \Bl_{\mm I}(X)$.
	The following theorem extends several results of Gaffney and Gassler \cite{GG} to arbitrary Noetherian local rings. 
	
	\begin{theorem}[Gaffney-Gassler intersection formulas]
		\label{thm:main}
		Assume \autoref{setup_general}.
		Then we have the following intersection product formulas and cutting results.
		\begin{enumerate}[\rm (i)]
			\item {\rm(Polar multiplicities):}
			$$
			m_i(I,R)
			\;=\;
			-\int_B \left(-E_I\right)^i \left(-E_\mm\right)^{d-i} \cdot \left[B\right]
			\quad \text{for $0\le i\le d-1$}.
			$$
			\item {\rm(Segre numbers):}
			$$
			c_i(I,R)
			\;=\;
			-\int_B \left(-E_I\right)^{i-1} \left(-E_\mm\right)^{d-i} \cdot \left[E_I\right]^{b^{-1}(U)} 
			\quad \text{ for $1\le i\le d-1$}.
			$$
			\item {\rm(Polar-Segre multiplicities):}
			$$
			\nu_i(I,R)
			\;=\;
			\int_B \left(-E_I\right)^{i-1} \left(-E_\mm\right)^{d-i} \cdot \left[E_I\right]^{b^{-1}(\{\mm\})} 
			\quad \text{ for $1\le i\le d-1$}.
			$$
			\item {\rm(Cutting with a general element of $I$):}
			Assume that $\kappa$ is infinite and $\dim(R/I)<d$.
			For a general element $f\in I$, we have
			$$
			m_i(I,R/fR) \;=\; m_{i+1}(I,R), \quad\;\;
			c_i(I,R/fR) \;=\; c_{i+1}(I,R), \quad\;\;
			\nu_i(I,R/fR) \;=\; \nu_{i+1}(I,R)
			$$
			for $1\leq i\leq d-2$.
			\item {\rm(Cutting with a general element of $\mm$):}
			Assume that $\kappa$ is infinite.
			For a general element $x\in\mm$, we have
			$$
			m_i(I,R/xR) \;=\; m_{i}(I,R), \quad\; c_i(I,R/xR) \;=\; c_{i}(I,R), \quad\; \nu_i(I,R/xR) \;=\; \nu_{i}(I,R)
			$$
			for $1 \le i \le d-2$.
		\end{enumerate}
	\end{theorem}
	
	\begin{proof}
		(i) Let $Y := \biProj\big(\gr_\mm(\Rees(I))\big) \subset B$.
		By \cite[\S 4]{KTBR}, \autoref{prop:basic} and \autoref{prop:joint-models}, we obtain
		\begin{align*}
			-\int_B \left(-E_I\right)^i \left(-E_\mm\right)^{d-i} \cdot \left[B\right] &\;=\; \int_B h_I^i h_\mm^{d-i-1} \cdot \left[Y\right] \\
			&\;=\; \int_Y c_1\left(\OO_Y(1,0)\right)^i c_1\left(\OO_Y(0,1)\right)^{d-i-1} \cdot \left[Y\right] \\
			&\;=\; m_i(I, R),
		\end{align*}
		where the last equality follows by \autoref{def_polar_segre}.
		
		(ii) Let $W := \biProj\big(\Rees(\mm \cdot\gr_I(R))\big) \subset B$ and $Z:=\biProj\big(\gr_\mm\left(\gr_I(R)\right)\big) \subset B$.
		From \cite[\S 4]{KTBR}, \autoref{prop:basic} and \autoref{prop:joint-models}, we get
		\begin{align*}
			-\int_B \left(-E_I\right)^{i-1} \left(-E_\mm\right)^{d-i} \cdot \left[E_I\right]^{b^{-1}(U)}  &\;=\;  -\int_B \left(-E_I\right)^{i-1} \left(-E_\mm\right)^{d-i} \cdot \left[W\right] \\
			&\;=\;  \int_B \left(-E_I\right)^{i-1} \left(-E_\mm\right)^{d-i-1} \cdot \left[Z\right] \\
			&\;=\; \int_B h_I^{i-1}h_\mm^{d-i-1} \cdot \left[Z\right] \\
			&\;=\; \int_Z c_1\left(\OO_Z(1,0)\right)^{i-1} c_1\left(\OO_Z(0,1)\right)^{d-i-1} \cdot \left[Z\right] \\
			&\;=\; c_i(I, R),
		\end{align*}
		where the last equality follows by \autoref{def_polar_segre}.
		
		(iii) For $1\leq i\leq d-1$, parts (i) and (ii) give
		\begin{align*}
			\nu_i(I, R) &\;=\; m_i(I, R) + c_i(I, R) \\
			&\;=\; \int_B \left(-E_I\right)^{i-1} \left(-E_\mm\right)^{d-i} \cdot \left(\left[E_I\right] - \left[E_I\right]^{b^{-1}(U)} \right) \\
			&
			\;=\; \int_B \left(-E_I\right)^{i-1} \left(-E_\mm\right)^{d-i} \cdot \left[E_I\right]^{b^{-1}(\{\mm\})}. 
		\end{align*}

		(iv) Let $\overline R:=R/fR$, which has dimension $d-1$.
		We have that $\OO_B(1,0)\cong \OO_B(-E_I) \cong I \OO_B$ is globally generated by sections induced from elements of $I$.
		Let $\sigma_f \in \HH^0(B, \OO_B(1,0))$ be the section induced by $f \in I$.
		Let $H := Z(\sigma_f)\subset B$ be the zero scheme determined by $\sigma_f$.
		By \autoref{lem:rees-superficial-specialization}, we have the isomorphism
		$$
		H \;\cong\; \biProj\big(T\big) \qquad \text{with} \qquad T \;:=\; \Rees_{\overline R}(I\overline R,\mm\overline R).
		$$
		The avoidance arguments of \autoref{lem:general-section-avoidance} and \autoref{lem_avoid_projection} imply that $H$ is an effective Cartier divisor on $B$,  that restricts to an effective Cartier divisor on $E_I$, and that it commutes with the decomposition of $[E_I]$ into its closed-fiber and punctured-spectrum parts:
		$$
		H \cdot \left[E_I\right]_{d-1}^{b^{-1}(\{\mm\})} \;=\; \left(H \cdot \left[E_I\right]_{d-1}\right)^{b^{-1}(\{\mm\})} \;=\;  \left(\left[\biProj\big(T/IT\big)\right]_{d-2}\right)^{b^{-1}(\{\mm\})} \;\in\; Z_{d-2}(B)
		$$
		and
		$$
		H \cdot \left[E_I\right]_{d-1}^{b^{-1}(U)} \;=\; \left(H \cdot \left[E_I\right]_{d-1}\right)^{b^{-1}(U)} \;=\;  \left(\left[\biProj\left(T/IT\right)\right]_{d-2}\right)^{b^{-1}(U)} \;\in\; Z_{d-2}(B).
		$$
		Therefore, part (i) yields 
		\begin{align*}
			m_{i+1}(I,R)
			\;&=\;
			-\int_B \left(-E_I\right)^{i+1} \left(-E_\mm\right)^{d-i-1} \cdot \left[B\right]\\
			\;&=\;
			-\int_B \left(-E_I\right)^{i} \left(-E_\mm\right)^{d-i-1} \cdot \left[\biProj(T)\right]\\
			\;&=\; m_i(I, R/fR).
		\end{align*}
		From part (ii), we obtain 
		\begin{align*}
			c_{i+1}(I,R)
			&\;=\;
			-\int_B \left(-E_I\right)^{i} \left(-E_\mm\right)^{d-i-1} \cdot \left[E_I\right]^{b^{-1}(U)} \\
			&\;=\;
			-\int_B \left(-E_I\right)^{i-1} \left(-E_\mm\right)^{d-i-1} \cdot \left[\biProj(T/IT)\right]^{b^{-1}(U)} \\
			&\;=\; c_i(I, R/fR).
		\end{align*}
		Since $\nu_i(I, R) = m_i(I, R) + c_i(I, R)$, the proof of this part is now complete. 
		
		(v) The proof of this part follows similarly to the one of part (iv).
	\end{proof}
	
	We recall the \emph{polar schemes} and \emph{Segre cycles} (see \cite{GG, CRPU}).
	Assume that $\kappa$ is infinite, and let $K_1,\ldots,K_d$ be general hyperplanes on $B_I$ corresponding to general elements $g_1,\ldots,g_d\in I$.
	Let
	$$
	P_i(I,X) \;:=\; b_I\left(K_1\cap\cdots\cap K_i\right)
	\qquad
	\text{and}
	\qquad
	\Lambda_i(I,X) \;:=\; (b_I)_*\left(\left[D_I \cap K_1\cap\cdots\cap K_{i-1}\right]_{d-i}\right),
	$$
	where the first image is scheme-theoretic.
	We now recover their multiplicity formulas for $1\leq i\leq d-1$.
	
	\begin{corollary}
		\label{cor:cycle-multiplicities}
		With notation as above, for all $1\leq i\leq d-1$, the following equalities hold:
		$$
		m_i(I,R) \;=\; e_{d-i}\left([P_i(I,X)]_{d-i}\right) \qquad \text{ and } \qquad
		c_i(I,R) \;=\; e_{d-i}\left(\Lambda_i(I,X)\right).
		$$
	\end{corollary}
	\begin{proof}
		Let $r:=d-i>0$ and $H_j:=(\widetilde{b}_\mm)^*(K_j)$.
		By \autoref{lem:general-section-avoidance}, we can assume that the following successive cuttings are all Cartier and yield the following effective cycles
		$$
		\Gamma_i \;:=\; H_1\cdots H_i\cdot[B]_d \;\in\; Z_r(B)
		\qquad \text{ and } \qquad
		\Sigma_i:=H_1\cdots H_{i-1}\cdot[E_I]_{d-1}^{b^{-1}(U)}
		\;\in\; Z_r(B).
		$$
		Successively applying the avoidance result of \autoref{lem_avoid_projection} implies that no component of $\Gamma_i$ is contained in $b^{-1}(V(I))$ and that no component of $\Sigma_i$ is contained in $b^{-1}(\{\mm\})$.
		
		The projection formula (at the level of cycles; see \cite[Proposition 1.4]{Fu}, \cite[\href{https://stacks.math.columbia.edu/tag/02RT}{Tag 02RT}]{stacks-project}) yields the equalities
		\begin{align*}
			b_*\left(\Sigma_i\right) \;=\; (b_I)_*\left((\widetilde{b}_\mm)_*\left(\Sigma_i\right)\right) 
			\;=\; (b_I)_*\left(K_1 \cdots K_{i-1} \cdot \left[D_I\right]_{d-1}^{b_I^{-1}(U)}\right) 
			\;=\; \Lambda_i(I, X),
		\end{align*}
		where the last equality follows because $r = d-i > 0$ and so $b_I$ pushes to zero a component of $[D_I]$ that is contained in $b_I^{-1}(\{\mm\})$.
		The natural isomorphisms $b^{-1}(X \setminus V(I)) \xrightarrow{\;\cong\;} b_I^{-1}(X \setminus V(I)) \xrightarrow{\;\cong\;} X \setminus V(I)$ yield the isomorphism
		$$
		\left(H_1 \cap \cdots \cap H_i\right)\mid_{b^{-1}(X \setminus V(I))} \;\cong\; \left(K_1 \cap \cdots \cap K_i\right)\mid_{b_I^{-1}(X \setminus V(I))}.
		$$
		Since no component of $H_1\cdots H_i\cdot[B]_d$ is contained in $b^{-1}(V(I))$ and no component of $K_1\cdots K_i \cdot \left[B_I\right]_d$ is contained in $b_I^{-1}(V(I))$, we obtain 
		$$
		b_*\left(\Gamma_i\right) \;=\; [P_i(I,X)]_r.
		$$
		Finally, \autoref{thm:main} and
		\autoref{lem:cycle-multiplicity-transfer} give
		$$
		m_i(I,R)
		\;=\;-\int_B(-E_\mm)^r\cdot\Gamma_i
		\;=\;e_r\left(b_*(\Gamma_i)\right)
		\;=\;e_r\bigl([P_i(I,X)]_r\bigr)
		$$
		and 
		$$
		c_i(I,R)
		\;=\;-\int_B(-E_\mm)^r\cdot\Sigma_i
		\;=\;e_r\left(b_*(\Sigma_i)\right)
		\;=\;e_r\bigl(\Lambda_i(I,X)\bigr).
		$$
		So the proof of the corollary is complete. 
	\end{proof}
	
	The following technical lemma justifies the multiplicity computation used above.
	\begin{lemma}
		\label{lem:cycle-multiplicity-transfer}
		Let $r>0$ and $\gamma\in Z_r(B)$ be a cycle with no component contained in $b^{-1}(\{\mm\})$. 
		Then
		$$
		e_r\left(b_*(\gamma)\right)\;=\; -\int_B (-E_\mm)^r \cdot\gamma.
		$$
	\end{lemma}
	\begin{proof}
		By linearity, we may assume $\gamma=[W]$, where $W\subset B$ is an integral closed subscheme of dimension $r$ such that  $W \not\subset b^{-1}(\{\mm\})$. 
		Let $V := b(W)$ and consider the blow-up
		$$
		\pi\colon \; \widetilde{V} :=\Bl_{\mm\OO_V}(V) \; \longrightarrow\; V
		$$
		with exceptional divisor $F$ determined by $\mm \OO_{\widetilde{V}} = \OO_{\widetilde{V}}(-F)$.
		Since
		$\mm\OO_W=\OO_W(-E_\mm|_W)$ is invertible, the universal property
		of the blow-up gives a proper factorization
		$
		b|_W\colon W \xrightarrow{f }\widetilde{V} \xrightarrow{ \pi } V
		$
		such that $f^*(F) = E_\mm\mid_W$.
		The projection formula and \cite[Theorem~3.1]{CutMon} yield
		$$
		-\int_B (-E_\mm)^r \cdot [W] \;=\; -\int_{\widetilde{V}} (-F)^r \cdot f_*\left([W]\right) \;=\; e_r\left(b_*\left([W]\right)\right),
		$$
		thus completing the proof of the lemma. 
	\end{proof}
	
	\subsection{Integral dependence under general hypersurface sections}
	\label{subsec:integral-dependence-general-sections}
	
	By \cite[Theorem~4.2]{PTUV}, Segre numbers characterize integral dependence in equidimensional and universally catenary Noetherian local rings.
	By \autoref{thm:main}(v), when the residue field is infinite, a general element $x\in\mm$ preserves the Segre numbers $c_i(I,R)$ for $1\leq i\leq d-2$.
	Together, these facts motivate the study of integral dependence under general hypersurface sections.
	
	For ideals $I\subseteq J$ in $R$, we consider the \emph{non-integral locus}
	$$
	N(I,J) \;:=\; \big\lbrace\pp\in \Spec(R) \,\mid\, \overline{I_\pp}\neq\overline{J_\pp} \big\rbrace \;\subset\; X.
	$$
	The following theorem, valid over an arbitrary Noetherian local ring with infinite residue field, shows that if $J(R/xR)$ is integral over $I(R/xR)$ for a general $x\in\mm$, then $\dim(N(I,J))\leq1$.
	
	\begin{theorem}
		\label{thm_gen_cut_mm}
		Assume \autoref{setup_general} with $\kappa$ infinite and $I \subseteq J$ two ideals in $R$. Suppose there exists a prime  $\mathfrak{p} \in X = \Spec(R)$ such that $\dim(R/\mathfrak{p}) \geq 2$ and $J_\mathfrak{p}$ is not integral over $I_\mathfrak{p}$. 
		Then, for a general element $x \in \mathfrak{m}$, the ideal $\overline{J} = J\overline{R}$ is not integral over $\overline{I} = I\overline{R}$ in the quotient ring $\overline{R} = R/xR$.
	\end{theorem}
	\begin{proof}
		Let $q : Y = \text{Proj}\left(\bigoplus_{n \geq 0} J^n\right) \rightarrow \text{Spec}(R)$ be the blow-up along the ideal $J$ and $E$ denote the exceptional divisor of $Y$. 
		Let $E_1,\ldots, E_k$ be the irreducible components of $E$.
		Let $\mathfrak{p}_i \in X = \Spec(R)$ denote the prime ideal corresponding to the generic point of the image $q(E_i)$.
		
		Since $I\mathscr{O}_Y \subseteq J\mathscr{O}_Y$ and $J\mathscr{O}_Y$ is an invertible sheaf, we can globally twist by its inverse to find a uniquely determined ideal sheaf $\mathfrak{b} \subseteq \mathscr{O}_Y$ such that
		$$
		I\mathscr{O}_Y \;=\;  \mathfrak{b} \cdot J\mathscr{O}_Y.
		$$
		The ideal sheaf $\mathfrak{b}$ describes the exact obstruction to the integrality of $J$ over $I$ (see, e.g., \cite[Corollary~1.2.5, \S8.2]{SwHu}).

		Let $W = V(\mathfrak{b}) \subset Y$ be the closed subscheme determined by $\mathfrak{b}$.
		Since $q(W) \subset X$ is closed and $\pp \in q(W)$, we obtain the inequality
		$$
		\dim(q(W)) \;\geq\; \dim(V(\mathfrak{p})) \;=\; \dim(R/\mathfrak{p}) \;\geq\; 2.
		$$
		
		Let
		$
		\mathfrak{S} :=  \lbrace \mathfrak{p}_i \mid 1 \le i \le k \text{ and } \mathfrak{p}_i \subsetneq \mathfrak{m} \rbrace
		$
		be the set of primes strictly contained in $\mathfrak{m}$ given as the image of the generic points of the irreducible components $E_i$.
		By prime avoidance, we obtain that 
		$$
		\dim\left(V(x) \cap q(W)\right) \;=\; \dim(q(W)) - 1 \;\geq\; 1
		$$
		and that 
		$$
		\text{$x \notin \mathfrak{p}_i$ \quad for any  \quad $\mathfrak{p}_i \in \mathfrak{S}$.}
		$$
		Thus we can fix a prime ideal $\mathfrak{q} \subsetneq \mm$ such that $\mathfrak{q} \in V(x) \cap q(W)$. 
		Since $\mathfrak{q} \in q(W)$, we can choose a point $w \in W$ in the blow-up $Y$ such that $q(w) = \mathfrak{q}$.
		Let $A := \mathscr{O}_{Y,w}$ be the local ring at $w$ and $\mathfrak{m}_w \subset \mathscr{O}_{Y,w}$ be its maximal ideal. 
		Since $x \in \mathfrak{q} = q(w)$, the pullback of $x$ into $A$ satisfies $x \in \mathfrak{m}_w$.
		
		By contradiction, assume that $\overline{J} = J \overline{R}$ is integral over $\overline{I} = I \overline{R}$ in the quotient ring $\overline{R}=R/xR$. 
		Considering the natural map $R/xR \rightarrow A/xA$ implies that $J(A/xA)$ is integral over $I(A/xA)$.
		
		Since $J\mathscr{O}_Y$ is an invertible sheaf, locally at $A$ it is generated by a single regular element $f \in A$, so $J A = (f)$. 
		The global equality $	I\mathscr{O}_Y =  \mathfrak{b} \cdot J\mathscr{O}_Y$ at the point $w$ turns into $I A = \mathfrak{b}_w f$. 
		Moreover, since $w$ was chosen to lie in $W = V(\mathfrak{b})$, we obtain $\mathfrak{b}_w \subseteq \mathfrak{m}_w$.
		Since $J(A/xA)$ is assumed to be integral over $I(A/xA)$, we have 
		$$
		f^n A \;\subseteq\; (\mathfrak{b}_w f) f^{n-1} A + xA \;=\; \mathfrak{b}_w f^n + xA.
		$$
		for some $n \ge 1$.
		This  implies $f^n = b f^n + ax $ for some $b \in \mathfrak{b}_w$ and $a \in A$. 
		Since $b \in \mathfrak{m}_w$, we obtain
		$$
		f^n \;=\; (1-b)^{-1}ax \;\in\; x A.
		$$
		We now divide the proof into two cases.

		{\sc Case 1:}
		Assume that $w \notin E$.
		Then $J A = A$  and so $A = (f^n) = xA$.
		This leads to the contradiction that $x$ is a unit in $A$ and $x \in \mm_w$.
		
		{\sc Case 2:}
		Assume that $w \in E$.
		Let $\mathfrak{Q} \in \Spec(A)$ be a minimal prime ideal lying over $xA$.
		The containment $f^nA \subset xA$ implies that $\mathfrak{Q}$ represents the generic point of some irreducible component $E_j$.
		Moving back to the base ring $R$, this means $x$ belongs to the prime ideal $\mathfrak{p}_j$. 
		However, since $w$ belongs to $E_j$, we should have $\mathfrak{p}_j \subseteq q(w) = \mathfrak{q}$.
		Finally, this leads to the contradiction
		$$
		x \;\in\; \pp_j \;\subseteq \; q(w) \;=\; \qqq \;\subsetneq\; \mm;
		$$
		indeed, we would  that $x$ belongs to the prime $\pp_j$ that lies in the set $\mathfrak{S}$.
		
		This concludes the proof of the theorem since we got a contradiction in both cases. 
	\end{proof}

	The next example shows that \autoref{thm_gen_cut_mm} is sharp.
	
	\begin{example}
		The assumption $\dim(R/\pp)\geq2$ in \autoref{thm_gen_cut_mm} is necessary.
		Let $\kk$ be an infinite field, $R:=\kk[x,y]_{(x,y)}$, $\mm:=(x,y)$, $
		I=(x^2,xy) \subset J:=(x,y)^2$ and $\pp:=(x)$.
		Then 
		$$
		\dim(R/\pp) \;=\;1, \qquad 
		I_\pp \;=\; xR_\pp\qquad\text{and}\qquad J_\pp \;=\; R_\pp.
		$$
		Thus $J_\pp$ is not integral over $I_\pp$.
		On the other hand, for a general element $\ell\in\mm$, we have the equality
		$$
		I \overline{R} \;=\; x^2\overline{R}
		\;=\; J\overline{R}
		$$
		in $\overline{R} = R/\ell R$.
	\end{example}
	
	\section{Integral dependence and Specialization}
	\label{sec:integral-dependence}
	
	In this section, we give a different proof of the following result that was proved in \cite[Theorem 3.3]{PTUV}.
	The statement of this result is a Principle of Specialization of Integral Dependence (PSID) type of argument \cite{GG}.
	A direct corollary of this result is a proof of the fact that Segre numbers detect integral dependence (see \cite[Theorem 4.2]{PTUV}).

	Throughout this section, we fix the following data. 
	
	\begin{setup}
		\label{setup_int_dep}
		Let $(R,\mm)$ be an equidimensional and universally catenary Noetherian local ring
		of dimension $d\geq2$. 
		Let $I\subset R$ be an ideal with height $\HT(I) \ge 1$.
		Let $t\in\mm$ be an element such that $\HT(t) = 1$ and $\HT(I, t) \ge 2$.
	\end{setup}

	\begin{theorem}[{\cite[Theorem 3.3]{PTUV}}]\label{thm:avoidance_main}
		Assume \autoref{setup_int_dep} and  that
		$$
		c_i(I,R/tR) \;\le\; c_i(I,R)\qquad\text{for all \;}1\leq i\leq d-1 .
		$$
		Then $t$, viewed in degree zero in $\gr_I(R)_0=R/I$, belongs to no minimal
		prime of $\gr_I(R)$.
	\end{theorem}
	
	In fact, we show the following more explicit extension of the above theorem. 
	It shows that Segre numbers are sensitive to the dimension of the projections of the \emph{bad components} (i.e., those lying over $V(t)$) of the exceptional divisor. 
	
	\begin{theorem}\label{thm:avoidance_main_explicit}
		Assume \autoref{setup_int_dep} and that there is a minimal prime of $\gr_I(R)$ that contains $t$. 
		Let $s$ be the maximal dimension
		$$
		s \;:=\; \max\big\lbrace \dim(R/Q \cap R) \mid \text{$Q$ minimal prime of $\gr_I(R)$ with $t \in Q$}\big\rbrace
		$$
		of the contractions to $R$ of the minimal primes of $\gr_I(R)$ containing $t$.
		Let $j := d-1-s$.
		Then we have the strict inequality
		$$
		c_j(I, R/tR) \;>\; c_j(I, R).
		$$
	\end{theorem}
	
	A direct consequence of the above theorem is the following algebraic version of the Principle of Specialization of Integral Dependence (PSID).
	
	\begin{corollary}\label{cor_PSID}
		Assume \autoref{setup_int_dep}.
		If $c_i(I, R)=c_i(I, R/tR)$ for all $1 \le i \le k$, then 
		$$
		\text{ for all \;$\pp \in \Ass_R\left(R/\overline{\,I\,}\right) \cap V(t),$\quad  we have $\HT(\pp)\ge k+2$.} 
		$$
		In particular, if $c_i(I, R)=c_i(I, R/tR)$ for all $1 \le i \le d-1$, then 
		$$
		f \in \overline{IR_t} \qquad 
		\text{ if and only if } \qquad f \in \overline{\,I\,}
		$$
		for all $f \in R$.
	\end{corollary}
	\begin{proof}
		By \cite[3.9 and 4.1]{MCADAM}, the associated primes of $R/\overline{\,I\,}$ come from contractions of minimal primes of $\gr_I(R)$ (see also \cite[\S9.6]{LAZ_POS_II}).
		So the claims follow from \autoref{thm:avoidance_main_explicit}.
	\end{proof}
	
	We now proceed to prove \autoref{thm:avoidance_main_explicit}.
	Let
	$$
	X \;:=\; \Spec(R),\qquad Z \;:=\; \Spec(R/I) \qquad \text{and} \qquad T \;:=\; \Spec(R/tR).
	$$
	Let
	$$
	\pi \colon \; B:=\Bl_Z(X) \;=\; \Proj(\Rees(I)) \;\longrightarrow\; X
	$$
	be the blow-up of $X=\Spec(R)$ along $Z$, and $E := E_ZX = \Proj(\gr_I(R))$ be the exceptional divisor. 	
	Throughout this section we shall use the fact that $\gr_I(R)$ is also equidimensional of dimension $d$.
	By utilizing the faithfully flat extension $R \rightarrow R[z]_{\mm R[z]}$, we may assume that the residue field of $R$ is not an algebraic extension of a finite field (see \autoref{subsec:infinite-residue-extension}).
	Let $\underline{H} = H_1,\ldots, H_d$ be a sequence of general hyperplanes on the blow-up $B=\Bl_Z(X)$, and denote by $\underline{g}=g_1,\ldots,g_d$ the associated sequence of elements in $I$ .

	By \cite[Theorem C]{CRPU}, the $i$-th polar scheme is given by 
	$$
	P_i(I, X) \; \;=\; \pi\left(H_1 \cap \ldots \cap H_i\right) \;=\; \Spec\left(R/(g_1,\ldots,g_i):_RI^\infty\right)
	$$
	and the $i$-th Segre cycle is given by
	$$
	\Lambda_i(I, X)  \;=\;  \pi_*\left(\big[E \cap H_1 \cap \cdots \cap H_{i-1}\big]_{d-i}\right) \;\in\; Z_{d-i}(X). 
	$$
	Moreover, we have that 
	$$
	c_i(I, R) \;=\; e_{d-i}\left(\Lambda_i(I, X)\right) \quad \text{ and } \quad \Lambda_i(I, X) = \Lambda_1\left(I, P_{i-1}(I, X)\right).
	$$
	
	We start with a simple specialization argument with respect to $T$.
	
	\begin{lemma}
		\label{lem_specialization_Seg}
		Assume that no irreducible components of $P_{i-1}(I, X)$ is contained in $T$. 
		Then we have the equality of cycles
		$$
		\Lambda_i\left(I, T\right) \;=\; \Lambda_1\left(I, \Spec(S/tS)\right),
		$$
		where $S=R/(g_1,\ldots,g_{i-1}):_RI^\infty$.
	\end{lemma}
	\begin{proof}
		Notice that $\Lambda_i(I, T) = \Lambda_1(I, P_{i-1}(I, T))$ and $\Lambda_1\left(I, \Spec(S/tS)\right) = \Lambda_1\left(I, \Spec(S/tS:_SI^\infty)\right) $.
		Since we have the equalities
		\begin{align*}
			P_{i-1}(I, T)  &\;=\; \Spec\big(R/(t, g_1,\ldots,g_{i-1}):_RI^\infty\big) \\
			& \;=\; \Spec\big(R/\left(( g_1,\ldots,g_{i-1}):_RI^\infty, t\right) :_R I^\infty\big) \;=\; \Spec(S/tS:_SI^\infty),
		\end{align*}
		the result of the lemma follows. 
	\end{proof}

	One of our main technical tools is the following explicit excess type of argument in the two dimensional case. 
	After cutting with suitably general hyperplanes on the blow-up, the general case will be deduced locally to this one. 
	
	\begin{proposition}[Surface excess]\label{prop_surface_excess}
		With notation as above, we assume that $\dim(R)=2$, $(0:_RI^\infty) = 0$, $\HT(t)=1$ and $(I,t)$ is $\mm$-primary.
		Then
		$$
		e\left(I,R/tR\right)
		\;\ge\;
		\sum_{\qq \in V(I),\; \HT(\qq)=1}
		e\left(I,R_\qq\right)e\left(t,R/\qq\right).
		$$
		If $E=E_ZX$ has an irreducible component mapping to the closed point
		$\mm$, then the inequality is strict.
	\end{proposition}
	\begin{proof}
		We may assume that $\HT(I) = 1$ because the statement of the proposition holds trivially when $I$ is $\mm$-primary.
		
		Choose a general element $g\in I$, and let $H \subset B$ be the corresponding hyperplane on the blow-up.  
		We choose $g$ so that: $g$ is a nonzerodivisor on $R$, $(g,t)$ is $\mm$-primary, $g$ is a reduction of $I(R/tR)$, and $H \subset B$ contains no irreducible component of $E$.
		
		The first Segre cycle is given explicitly by 
		$$
		\Lambda_1(I, X) \;=\; \pi_*\left(\left[E\right]\right) \;=\; \sum_{\qq\in V(I),\; \HT(\qq)=1}
		e\left(I,R_\qq\right)\left[V(\qq)\right]
		$$
		(see \cite[Theorem~C(ii)]{CRPU}, \cite[Proposition~2.4(a)]{PTUV} and \cite[Proposition~2.3(iii)]{AM}).
		Let $D = V(g)$ be the effective Cartier divisor on $X$ corresponding to $g$.
		We have the equality of Cartier divisors 
		$$
		\pi^*(D) \;=\; E + H
		$$
		on $B$ (see \cite[Proposition 5.8]{CRPU}, \cite[Lemma 2.4]{Fu}).
		Then the projection formula yields the equality of cycles
		$$
		\left[D\right]\;=\; \pi_*\left(\left[\pi^*(D)\right]\right) \;=\; \pi_*\left(\left[E\right]\right) + \pi_*\left(\left[H\right]\right) \;\in\; Z_1(X)
		$$
		on $X$ (see \cite[Proposition 1.4]{Fu}, \cite[\href{https://stacks.math.columbia.edu/tag/02RT}{Tag 02RT}]{stacks-project}).
		Since $(g,t)$ is $\mm$-primary, no component of $D$ is contained in $T$.
		By \cite[Corollary 1.2.14]{FOV}, we get
		$$
		e\left((g,t),R\right) \;=\; e\left(t, R/gR\right) \;=\; e\left(t, \left[D\right]\right)
		$$
		and 
		$$
		e\left(g, R/tR\right) \;\ge\; e\left((g,t),R\right).
		$$
		Thus, since $g$ is a reduction of
		$I(R/tR)$, we obtain
		$$
		e(I,R/tR) \;=\; e\left(g, R/tR\right) \;\geq \; e\left(t, \left[D\right]\right).
		$$
		On the other hand, we have
		$$
		e\left(t, \pi_*\left(\left[E\right]\right)\right)
		\;=\;
		\sum_{\qq\in V(I),\; \HT
			(\qq)=1}
		e\left(I,R_\qq\right)e\left(t,R/\qq\right).
		$$
		By combining the above equations, we obtain 
		\begin{align}
			\begin{split}
				\label{eq_surface_excess}
				e(I,R/tR) &\;\ge\; e\left(t, \left[D\right]\right) \\
				&\;=\; e\left(t, \pi_*\left(\left[E\right]\right)\right) + e\left(t, \pi_*\left(\left[H\right]\right)\right) \\
				&\;=\; \sum_{\qq\in V(I),\; \HT
					(\qq)=1}
				e\left(I,R_\qq\right)e\left(t,R/\qq\right) \;+\; e\left(t, \pi_*\left(\left[H\right]\right)\right).
			\end{split}
		\end{align}
		This proves the weak inequality.
		
		Now assume that $E$ has an irreducible component $F$ over the closed point $\mm$.
		We show that this assumption implies the nonvanishing $H\neq \varnothing$.
		Since $\OO_F(1) = \OO_B(1) \otimes_{\OO_B} \OO_F$ is ample over $\Spec(R/\mm)$, it has positive degree. 
		By assumption $H$ restricts to a general section of $\OO_F(1)$, and so we get $F \cap H \neq \varnothing$.
		In particular, we obtain $H\neq \varnothing$.
		
		Let $C \subset H$ be an irreducible component of $H$.
		It necessarily follows that $C$ is not contained $E = B \times_X \Spec(R/I)$ because by assumption $H$ contains no irreducible component of $E$. 
		This implies that $C$ is not contained in the special fiber 
		$$
		B \times_{X} \Spec(R/\mm) \;\subseteq\; B \times_{X} \Spec(R/I) = E.
		$$
		Therefore, $\pi_*\left(\left[C\right]\right) \neq 0$, and so $\pi_*\left(\left[H\right]\right) \neq 0$. 
		Since no component of $\pi_*
		\left(\left[H\right]\right)$ is contained in $T$, we get
		$e\left(t, \pi_*\left(\left[H\right]\right)\right) > 0$.
		Finally, by \autoref{eq_surface_excess}, the strict inequality follows. 
	\end{proof}
	
	The above proposition already yields the result of \autoref{thm:avoidance_main_explicit} in the following case. 
	
	\begin{corollary}
		\label{cor_main_c_1}
		With notation as above, we assume $(0:_RI^\infty)=0$, $\HT(t)=1$ and $\HT(I,t) \ge 2$. 
		Suppose that $E$ has an irreducible component $F$ such that
		$$
		\pi(F) \subset T,\qquad \pp_0 \in \Spec(R) \text{ with } V(\pp_0) = \pi(F),\qquad \HT(\pp_0)=2.
		$$
		Then  $c_1(I, R/tR) > c_1(I, R)$.
	\end{corollary}
	\begin{proof}
		We consider the following Segre cycles
		$$
		\Lambda_1(I, X) \;=\; \sum_{\qq\in V(I), \; \HT(\qq)=1}
		e\left(I,R_\qq\right)\left[V(\qq)\right]
		$$
		and 
		$$
		\Lambda_1(I, T)  \;=\; \sum_{\pp\in V(I,t), \; \HT(\pp)=2}
		e\left(I,R_\pp/tR_\pp\right)\left[V(\pp)\right].
		$$
		Recall that $c_1(I, R) = e\left(\Lambda_1(I, X)\right)$ and $c_1(I, R/tR) = e\left(\Lambda_1(I, T)\right)$.
		Since $\HT(I,t) \ge 2$, no component of $\Lambda_1(I, X)$ is contained in $T$.
		Hence we can consider the cycle
		\begin{align*}
			T \cdot \Lambda_1(I, X) &\;=\; \sum_{\qq\in V(I),\; \HT(\qq)=1}
			e\left(I,R_\qq\right) \sum_{\pp\in V(I,t),\; \HT(\pp)=2}
			e\left(t,R_\pp/\qq R_\pp\right)\left[V(\pp)\right] \\
			&\;=\; \sum_{\pp\in V(I,t),\; \HT(\pp)=2}
			\left(
			\sum_{\qq\in V(I),\; \qq\subset \pp,\; \HT(\qq)=1} e\left(I,R_\qq\right)
			e\left(t,R_\pp/\qq R_\pp\right)
			\right)
			\left[V(\pp)\right].
		\end{align*}
		Therefore, by applying \autoref{prop_surface_excess} to each localization $R_\pp$ with $\pp \in V(I,t)$ and $\HT(\pp) = 2$, we obtain
		\begin{align*}
			{\rm Coeff}_{[V(\pp)]}\left(\Lambda_1(I, T)\right) &\;=\; e\left(I,R_\pp/tR_\pp\right)\\ &\;\ge\; \sum_{\qq\in V(I),\; \qq\subset \pp,\; \HT(\qq)=1} e\left(I,R_\qq\right)
			e\left(t,R_\pp/\qq R_\pp\right)\\
			& \;=\; {\rm Coeff}_{[V(\pp)]}\left(T\cdot\Lambda_1(I, X)\right).
		\end{align*}
		By assumption, the exceptional divisor $E_{\pp_0} = \Proj\left(\gr_I(R_{\pp_0})\right)$ has an irreducible component lying over the closed point of $\Spec(R_{\pp_0})$, and so \autoref{prop_surface_excess} yields the strict inequality 
		$$
		{\rm Coeff}_{[V(\pp_0)]}\left(\Lambda_1(I, T)\right) \;>\; {\rm Coeff}_{[V(\pp_0)]}\left(T\cdot\Lambda_1(I, X)\right).
		$$
		Thus $e\left(\Lambda_1(I, T)\right) > e\left(T \cdot \Lambda_1(I, X)\right)$. 
		By combining this with \cite[Corollary 1.2.12]{FOV}, we get
		$$
		c_1(I, R/tR) \;=\; e\left(\Lambda_1(I, T)\right) \;>\; e\left(T \cdot \Lambda_1(I, X)\right) \;\ge\;  e\left(\Lambda_1(I, X)\right) \;=\; c_1(I, R).
		$$
		So the proof of the corollary is complete.
	\end{proof}
	
	Finally, we are ready for the proof of \autoref{thm:avoidance_main_explicit}.
	
	\begin{proof}[Proof of \autoref{thm:avoidance_main_explicit}]
		Let $Q$ be a minimal prime of $\gr_I(R)$ containing $t$ such that the contraction $\pp_0 := Q \cap R \in \Spec(R)$ has the maximal  possible dimension $\dim(R/\pp_0) = s$.
		Note that $Q$ corresponds to an irreducible component $F$ of $E=\Proj(\gr_I(R))$ because $\dim(\gr_I(R)/\gr_I(R)_+) = \dim(R/I) < \dim(\gr_I(R)) = d$.
		Let $j := d-1-s$.
		Since $\HT(I, t) \ge 2$, we have $0 \le s \le d-2$ and so $1 \le j \le d-1$.
		
		We analyze the $j$-th Segre cycle
		$$
		\Lambda_j(I, X) \;=\; \Lambda_1\left(I, Y\right),
		$$
		where $Y := P_{j-1}(I, X) = \pi\left(H_1 \cap \cdots \cap H_{j-1}\right) = \Spec(S)$ and $S := R/(g_1,\ldots,g_{j-1}):_RI^\infty $.
		
		We now check the necessary conditions to apply \autoref{cor_main_c_1} to $S$:
		\begin{enumerate}[(a)]
			\item	Let $\kappa(\pp_0):=R_{\pp_0}/\pp_0R_{\pp_0}$.
			The fact that $\gr_I(R)$ has a minimal prime contracting to $\pp_0$ is equivalent to saying that $I_{\pp_0}$ has maximal analytic spread $\ell(I_{\pp_0}):=\dim(\gr_I(R) \otimes_R \kappa(\pp_0)) = \dim(R_{\pp_0})=d-s$. 
			Since we may assume that $g_1,\ldots,g_d$ are general elements in $I_{\pp_0}$ (see \cite[Lemma 2.6]{PTUV}), it follows that 
			$$
			\dim\left(\gr_I(S) \otimes_R \kappa(\pp_0)\right) \;=\; d-s-j+1=2 \quad \text{ and } \quad \dim(S_{\pp_0}) \;=\; d-s-j+1=2
			$$
			(see \cite[\S 8.5]{SwHu} and \cite[Remark 5.11]{CRPU}).
			This implies that $\pp_0 \supset (g_1,\ldots,g_{j-1}):_RI^\infty$, that $\HT(\pp_0 S)=2$ and that there is a minimal prime of $\gr_I(S)$ lying over $\pp_0S$.
			We have that $S$ is equidimensional of dimension $\dim(S) = d - j + 1 = s+2$ (see \cite[Remark 5.11]{CRPU}).
			Recall also that $t \in \pp_0$.
			
			\item By applying \autoref{lem_avoid_projection} successively to the irreducible components of the blow-up $B$, we obtain that no component of $Y = \pi(H_1 \cap \cdots \cap H_{j-1})$ is contained in $T$, and so $\HT(tS) \ge 1$.
			
			\item 	The components of $\Lambda_j(I, X) = \pi_*([E \cap H_1 \cap \cdots \cap H_{j-1}]_{d-j})$ have dimension $d-j=s+1$.
			Let $C \subset E$ be an irreducible component of $E$.
			If $\pi(C)$ is not contained in $T$, then by applying \autoref{lem_avoid_projection} successively to $C$, we obtain that no component of $ \pi(C \cap H_1 \cap \cdots \cap H_{j-1})$ is contained in $T$. 
			If $\pi(C)$ is contained in $T$, then $\dim(\pi(C \cap H_1 \cap \cdots \cap H_{j-1})) \le \dim(\pi(C)) \le s$.
			By definition of pushforward, this implies that no component of $\Lambda_1\left(I, Y\right)=\Lambda_j(I, X)$ is contained in $T$.
			Since 	
			$$
			\Lambda_1(I, Y) \;=\; \sum_{\qq\in V(IS), \; \HT(\qq)=1}
			e\left(I,S_\qq\right)\left[V(\qq)\right],
			$$
			it follows that $\HT(IS, tS) \ge 2$.
		\end{enumerate}
		From \autoref{cor_main_c_1}, we obtain
		$$
		c_1(I, S/tS) \;>\; c_1(I, S).
		$$
		Finally, the equalities of cycles $\Lambda_1\left(I, \Spec(S/tS)\right) = \Lambda_j(I, T)$ (see \autoref{lem_specialization_Seg}) and $\Lambda_1(I,  Y) = \Lambda_j(I, X)$ yield the strict inequality 
		$$
		c_j(I, R/tR) \;=\; c_1(I, S/tS) \;>\; c_1(I, S) \;=\; c_j(I, R).
		$$
		So the proof of the theorem is complete. 
	\end{proof}

	\section*{Acknowledgments}
	
	YCR received support from NSF grant DMS-2502321 and Simons Foundation Travel Support for
	Mathematicians Award MPS-TSM-00013551.
	SDC received support from NSF grant DMS-2348849.

	\bibliographystyle{amsalpha}

\begin{thebibliography}{PTUV20}
		
		\bibitem[AM97]{AM}
		R{\"u}diger Achilles and Mirella Manaresi,
		\emph{\href{https://doi.org/10.1007/s002080050128}{Multiplicities of a
				bigraded ring and intersection theory}}, Math. Ann. \textbf{309} (1997),
		no.~4, 573--591.
		
		\bibitem[AR01]{AR}
		R{\"u}diger Achilles and S{\l}awomir Rams,
		\emph{\href{https://doi.org/10.1007/PL00000509}{Intersection numbers, Segre
				numbers and generalized Samuel multiplicities}}, Arch. Math. (Basel)
		\textbf{77} (2001), no.~5, 391--398.
		
		\bibitem[CR24]{CidPolar}
		Yairon Cid-Ruiz,
		\emph{\href{https://doi.org/10.1093/imrn/rnae163}{Polar multiplicities and
				integral dependence}}, Int. Math. Res. Not. IMRN \textbf{2024} (2024),
		no.~17, 12201--12218.
		
		\bibitem[CM25]{CutMon}
		Steven~Dale Cutkosky and Jonathan Monta{\~n}o,
		\emph{\href{https://doi.org/10.1112/jlms.70338}{Multiplicities and degree
				functions in local rings via intersection products}}, J. Lond. Math. Soc.
		\textbf{112} (2025), no.~5, e70338.
		
		\bibitem[CM26]{CutMon2}
		\bysame, \emph{\href{https://doi.org/10.1007/s13348-026-00519-w}{Degree
				functions of graded families of ideals}}, Collect. Math. (2026),
		published online.
		
		\bibitem[CRPU25]{CRPU}
		Yairon Cid-Ruiz, Claudia Polini, and Bernd Ulrich,
		\emph{\href{https://arxiv.org/abs/2405.07000}{Multidegrees, families, and
				integral dependence}}, preprint, 2025.
		
		\bibitem[Cut26]{Cut26}
		Steven~Dale Cutkosky,
		\emph{\href{https://arxiv.org/abs/2603.06844}{Multiplicities of graded
				families of ideals on Noetherian local rings}}, arXiv:2603.06844, 2026.
		
		\bibitem[FM01]{FleMan}
		Hubert Flenner and Mirella Manaresi,
		\emph{\href{https://doi.org/10.1007/PL00004900}{A numerical characterization
				of reduction ideals}}, Math. Z. \textbf{238} (2001), no.~1, 205--214.
		
		\bibitem[FOV99]{FOV}
		Hubert Flenner, Liam O'Carroll, and Wolfgang Vogel,
		\emph{\href{https://doi.org/10.1007/978-3-662-03817-8}{Joins and
				intersections}}, Springer Monographs in Mathematics, Springer, Berlin, 1999.
		
		\bibitem[Ful98]{Fu}
		William Fulton,
		\emph{\href{https://doi.org/10.1007/978-1-4612-1700-8}{Intersection Theory}},
		second ed., Ergebnisse der Mathematik und ihrer Grenzgebiete (3), vol.~2,
		Springer-Verlag, Berlin, 1998.
		
		\bibitem[GG99]{GG}
		Terence Gaffney and Robert Gassler,
		\emph{\href{https://arxiv.org/abs/alg-geom/9611002}{Segre numbers and
				hypersurface singularities}}, J. Algebraic Geom. \textbf{8} (1999), no.~4,
		695--736.
		
		\bibitem[KT94]{KTBR}
		Steven~L. Kleiman and Anders Thorup,
		\emph{\href{https://doi.org/10.1006/jabr.1994.1182}{A geometric theory of the
				Buchsbaum--Rim multiplicity}}, J. Algebra \textbf{167} (1994), no.~1,
		168--231.
		
		\bibitem[KT96]{KTMixed}
		\bysame, \emph{\href{https://arxiv.org/abs/alg-geom/9412019}{Mixed
				Buchsbaum--Rim multiplicities}}, Amer. J. Math. \textbf{118} (1996), no.~3,
		529--569.
		
		\bibitem[Laz04]{LAZ_POS_II}
		Robert Lazarsfeld,
		\emph{\href{https://doi.org/10.1007/978-3-642-18810-7}{Positivity in
				algebraic geometry II: Positivity for vector bundles, and multiplier
				ideals}}, Ergebnisse der Mathematik und ihrer Grenzgebiete (3), vol.~49,
		Springer-Verlag, Berlin, 2004.
		
		\bibitem[Mat86]{MatCRT}
		Hideyuki Matsumura,
		\emph{\href{https://doi.org/10.1017/CBO9781139171762}{Commutative ring
				theory}}, Cambridge Studies in Advanced Mathematics, vol.~8,
		Cambridge University Press, Cambridge, 1986, translated from the Japanese
		by M.~Reid.
		
		\bibitem[McA83]{MCADAM}
		Stephen McAdam,
		\emph{\href{https://doi.org/10.1007/BFb0071575}{Asymptotic prime divisors}},
		Lecture Notes in Math., vol.~1023, Springer-Verlag, Berlin, 1983.
		
		\bibitem[PTUV20]{PTUV}
		Claudia Polini, Ng{\^o}~Vi{\^e}t Trung, Bernd Ulrich, and Javid Validashti,
		\emph{\href{https://doi.org/10.1007/s00208-020-02059-5}{Multiplicity sequence
				and integral dependence}}, Math. Ann. \textbf{378} (2020), 951--969.
		
		\bibitem[Ree61]{Rees61}
		David Rees,
		\emph{\href{https://doi.org/10.1017/S0305004100034800}{$\mathfrak{a}$-transforms
				of local rings and a theorem on multiplicities of ideals}},
		Proc. Cambridge Philos. Soc. \textbf{57} (1961), no.~1, 8--17.
		
		\bibitem[Ree84]{Rees84}
		David Rees,
		\emph{\href{https://doi.org/10.1112/jlms/s2-29.3.397}{Generalizations of
				reductions and mixed multiplicities}}, J. London Math. Soc. (2)
		\textbf{29} (1984), no.~3, 397--414.
		
		\bibitem[SH06]{SwHu}
		Irena Swanson and Craig Huneke,
		\emph{\href{https://www.math.purdue.edu/~iswanso/book/index.html}{Integral
				closure of ideals, rings, and modules}}, London Mathematical Society Lecture
		Note Series, vol. 336, Cambridge Univ. Press, Cambridge, 2006.
		
		\bibitem[{Sta}]{stacks-project}
		The {Stacks project authors}, \emph{\href{https://stacks.math.columbia.edu}{The
				Stacks project}}, \url{https://stacks.math.columbia.edu}.
		
		\bibitem[Tei73]{TEISSIER_CYC}
		Bernard Teissier,
		\emph{\href{https://numdam.org/item/AST_1973__7-8__285_0/}{Cycles
				\'evanescents, sections planes et conditions de Whitney}},
		Singularit\'es \`a Carg\`ese, Ast\'erisque, nos.~7--8,
		Soc. Math. France, Paris, 1973, pp.~285--362.
		
		\bibitem[Tei80]{TEISSIER_RES2}
		\bysame, \emph{\href{https://doi.org/10.1007/BFb0085881}{R\'esolution
				simultan\'ee II: R\'esolution simultan\'ee et cycles \'evanescents}},
		S\'eminaire sur les singularit\'es des surfaces (M.~Demazure,
		H.~C.~Pinkham, and B.~Teissier, eds.), Lecture Notes in Math., vol.~777,
		Springer-Verlag, Berlin, 1980, pp.~82--146.
		
		\bibitem[Tho90]{Thorup}
		Anders Thorup, \emph{\href{https://doi.org/10.1007/BFb0084049}{Rational
				equivalence theory on arbitrary Noetherian schemes}}, Enumerative Geometry
		(Sitges, 1987), Lecture Notes in Math., vol. 1436, Springer, Berlin, 1990,
		pp.~256--297.
		
		\bibitem[TV07]{TV}
		Ng{\^o}~Vi{\^e}t Trung and Jugal~K. Verma,
		\emph{\href{https://doi.org/10.1090/S0002-9947-07-04054-8}{Mixed
				multiplicities of ideals versus mixed volumes of polytopes}},
		Trans. Amer. Math. Soc. \textbf{359} (2007), no.~10, 4711--4727.
		
	\end{thebibliography}

\end{document}